\documentstyle[12pt,amscd]{amsart}
\newtheorem{Theorem}{Theorem}[section]
\newtheorem{Lemma}[Theorem]{Lemma}
\newtheorem{Corollary}[Theorem]{Corollary}
\newtheorem{Proposition}[Theorem]{Proposition}

\newtheorem{Definition}[Theorem]{Definition}

\def\depth{\operatorname{depth}}
\def\reg{\operatorname{reg}}
\def\Ext{\operatorname{Ext}}
\def\geom{\operatorname{g-reg}}
\def\hdeg{\operatorname{hdeg}}
\def\Ann{\operatorname{Ann}}

\def\sk{\smallskip\par}
\def\mm{{\mathfrak m}}

\def\qq{{\mathfrak q}}
\def\RR{{\mathfrak R}}

\def\Mc{{\mathcal M}}

\def\FilM{{\mathbb M}}
\def\ZZ{{\mathbb Z}}
\def\NN{{\mathbb N}}

\begin{document}
\title{Castelnuovo-Mumford regularity of \\
 associated graded modules  and   fiber cones \\ of  filtered modules}
\thanks{Both authors were partially supported by NAFOSTED
(Vietnam)}
\maketitle
\begin{center} 
LE XUAN DUNG\\
Department of  Natural Sciences,  Hong Duc University\\
307 Le Lai, Thanh Hoa, Vietnam\\
E-mail: lxdung27@@gmail.com\\ [15pt] 
and\\[15pt] 
LE TUAN HOA \\ 
 Institute of Mathematics\\ 18 Hoang Quoc Viet, 10307 Ha Noi, Vietnam\\
 E-mail: lthoa@@math.ac.vn
  \end{center}
 
\begin{abstract}  We give bounds on the Castelnuovo-Mumford regularity of the
associated graded module of  an arbitrary good filtration  and of its fiber cone. These bounds  extend previous results of  Rossi-Trung-Valla and Linh.
\end{abstract}

\date{}
\section*{Introduction} \sk
Let $A$ be  a local ring with the maximal ideal $\mm$ and  $M$ be a finitely generated  $A$-module. Denote by $G_I(M) = \oplus_{n\geq 0}I^nM/I^{n+1}M $  the associated graded module of $M$ with respect to  an $\mm$-primary ideal $I$. It is well-known that one can use the Castelnuovo-Mumford regularity $\reg(G_I(M))$ of $G_I(M)$ to bound other invariants of $M$  such as  Hilbert  coefficients, the postulation number, ... (see, e.g., \cite{T}, \cite{Va}, \cite{RTV}). Therefore, bounding  $\reg(G_I(M))$ is an important problem. This problem was  completely solved by  Rossi, Trung and Valla  for the case $M=A$ and $I= \mm$  in \cite{RTV} and  by Linh  for the general case \cite{L1}, where bounds on $\reg(G_I(M))$ were given in terms of the dimension $d$ of $M$ and of the so-called extended degree $D(I,M)$ (see Definition \ref{A1}). 

Another important subject associated to  $I$ is the so-called  fiber cone $F_\mm(I) =  \oplus_{n\geq 0}I^n/\mm I^n$. Inspired by the mentioned results of \cite{RTV} and \cite{L1}, it is natural to ask whether one can bound the Castelnuovo-Mumford regularity  $\reg(F_\mm(I))$ of $F_\mm(I)$ in terms of $D(I,A)$ and $d$?  One of the main ideas of  \cite{RTV} and \cite{L1} is to reduce the problem of bounding  the Castelnuovo-Mumford regularity $\reg(G_I(M))$  in the case $\depth M>0$ to bounding the  so-called geometric  Castelnuovo-Mumford regularity of $G_I(M)$. However this approach does not work for fiber cones, because a similar result to \cite[Theorem 5.2]{H} does not hold in this case.  Fortunately, thank to a recent work by Rossi and Valla \cite{RV}, we can  use the associated graded modules of  some good filtrations of $M$ to solve the problem. 

The notion of good $I$-filtrations $\FilM = \{ M_n\}_{n\ge 0}$ of $M$ was already considered  in \cite{Bour} and \cite{AM} (see Definition \ref{Fil-Def}). Even in the case $M=A$ the class good filtrations of submodules is larger than the usual class of good filtrations of ideals. Recently, Rossi and Valla showed in \cite{RV} that one can not only extend many classical results to a filtered module (i.e. a module with a good filtration), but also use filtered modules to study the Hilbert function of fiber cones. Our work here once more shows the usefulness of this concept  in studying fiber cones. 

In order to bound $\reg(F_\mm(I))$ we first extend results of \cite{RTV,L1} to the associated graded module $G(\FilM) = \oplus_{n\ge 0}M_n/M_{n+1}$ of $\FilM$.  The techniques we use are similar to that  in  \cite{RTV,L1}. The new point here is that we also use the maximal generating degree   $r(\FilM)$ of  the graded module $G(\FilM)$ over the standard graded ring $G_I(A)$ (see Definition  \ref{A2}) and clarify its influence to the Castelnuovo-Mumford regularity of  associated graded modules by a hyperplane section or by passing to  $\bar{M} = M/H^0_\mm(M)$. As a consequence, we see that the bounds in this general case (see Theorem \ref{A3}) look similar to that in the classical case in \cite{L1}.  The bounds are given in terms of  $r(\FilM), \ d$ and in terms of  $D(I,M)$, which does not depend on the filtration $\FilM$. Then using an exact sequence connecting fiber cones and associated graded modules of some good filtrations  and a relation between Hilbert coefficients of these filtrations given in \cite{RV} we can also bound the Castelnuovo-Mumford regularity of a fiber cone  associated to any good filtration  in terms of $d,\ D(I,A)$ and $r(\FilM)$ (see Theorem \ref{F3}). 

In this paper we also consider the graded case.  For this situation we can apply the approach of \cite{RTV} and \cite{L1}  only in the case $I$ being generated by homogeneous elements of the same degree.  A new point here is that  we can give a bound on  $\reg(G(\FilM))$ in terms of  $\reg(M)$ (see Theorem \ref{C3}).  Note that $\reg(M)$ is in general not only much smaller  than an arbitrary extended degree $D(I,M)$ of $M$ (see \cite{DGV, Va, Na}), but  it is  also much easier to compute. In the general case the above mentioned  approach is not applicable, because $I$ does not contain a generic homogeneous element. To overcome the difficulty we pass to the localization and first give a bound in terms of the so-called  homological degree (Theorem \ref{B3}). Combining with a recent result of \cite{CHH} we can then get a bound in terms of $\reg(M)$ (see Theorem \ref{B5}),  provided $A$ is a polynomial ring over a field. However the bound in this general setting is much worse than the one in Theorem \ref{C3}.

The paper is divided into three sections. In Section \ref{Local} we bound  the Castelnuovo-Mumford regularity of a filtered module in the local case. The graded case is considered in Section \ref{Graded}. In the last section  \ref{Fiber} we derive bounds on  the Castelnuovo-Mumford regularity of the fiber cone of a filtered module.

\section{Regularity of associated graded rings: Local case} \label{Local}

Let $A$ be a Noetherian  local ring with an infinite residue field $K:= A/\mm$ and $M$ a finitely generated $A$-module. (Although the assumption $K$ being  infinite is not essential, because we can tensor $A$ with $K(t)$.)  First we recall some basic facts on filtered modules from \cite[Section III.3]{Bour}, \cite[Chapter 10]{AM} and \cite[Section 1]{RV}.

\begin{Definition} \label{Fil-Def} {\rm Given a proper ideal $I$. A chain of submodules
$$\FilM:\ M = M_0 \supseteq M_1 \supseteq M_2 \supseteq \cdots \supseteq M_n \supseteq  \cdots $$
 is called an {\it $I$-filtration} of $M$ if $IM_i \subseteq  M_{i+1}$ for all $i$, and  a {\it good $I$-filtration} if $IM_i =  M_{i+1}$ for all sufficiently large $i$. A module $M$ with a filtration is called a} filtered module.
 \end{Definition}
 
 Thus $\{I^nM\}$ is a good $I$-filtration. Note that the notion of $I$-filtration of submodules is different from that of filtration of ideals. The sequence $A \supset \mm \supset \mm I \supset \mm I^2 \supset  \cdots $ is a good $I$-filtration of submodules of $A$, but it is in general not a filtration of ideals in  $A$, because $\mm^2 \not\subseteq  \mm I$.  

If $N$ is a submodule of $M$, then by Artin-Rees Lemma, the sequence $\{N\cap M_n\}$ is a good $I$-filtration of $N$ and we will denote it by $\FilM \cap N$. The sequence $\{M_n+N/N\}$ is a good $I$-filtration of $M/N$ and will be denoted by $\FilM/N$.

In this paper we always assume that $I$ is an $\mm$-primary ideal and $\FilM$ is a good $I$-filtration. The {\it associated graded module}  to the filtration $\FilM$ is defined by
$$G(\FilM) = \bigoplus _{n\geq 0}M_n/M_{n+1} .$$
We also say that $G(\FilM)$ is the associated ring of the filtered module $M$. This is a finitely generated graded module over the standard graded ring $G:= G_I(A) := \oplus_{n\ge 0}I^n/I^{n+1}$ (see \cite[Proposition III.3.3]{Bour}). 

\begin{Definition}\label{A2}{\rm Let  $\FilM$ be a good $I$-filtration of $M$. We set
$$r = r(\FilM) = \min \{t \geq 0 \mid   M_{n+1} = IM_n \ \ \text{for all} \ \ n \geq  t \}.$$}
\end{Definition}

In particular, in the $I$-adic case, $r(\{I^nM\}) = 0$. Note that $r$ is always finite, and $M_{r+j} = I^jM_r$ for all $j\ge 0$. This means $\{ M_n\}_{n\ge r}$ is of form of an $I$-adic filtration of $M_r$. In other words, $r$ is the largest generating degree of $G(\FilM)$ as a graded module over $G$. Together with well-known facts on the Castelnuovo-Mumford regularity (see, e.g. \cite[Paragraph 1.3]{BS}) this explains why $r$ will naturally occur in our bounds on $\reg(G(\FilM))$, which is recalled below.

The {\it Castelnuovo-Mumford regularity} of a finitely generated graded module $E = \oplus _{n \in \ZZ} E_n$ over a standard $\NN$-graded ring $R = \oplus _{n\geq 0} R_n$  is defined as the number 
$$\reg(E) = \max \{ a_i(E) + i \mid  i \geq  0 \},$$
where 
$$ a_i(E) =
\begin{cases}
\sup \{n|\ H_{R_+}^i(E)_n \ne 0 \}  \ \ \ \ \text{if} \  \ H_{R_+}^i(E) \ne 0 ,\\
 -\infty \  \hskip 3.6cm \text{if} \ \ H_{R_+}^i(E) = 0,
\end{cases}
$$
 and $R_+ = \oplus _{n>  0} R_n.$ We will write $\reg(G(\FilM))$ to mean the Castelnuovo-Mumford regularity of $G(\FilM)$ being a graded module over the standard graded ring $G$.
 
 An element $x\in R_1$ is called  (linear) {\it filter-regular } on $E$ if $(0:_Ex)_n = 0$ for all $n\gg 0$.  From this one can show that $0:_Ex \subseteq H^0_{R_+}(E)$, and hence  $(0:_Ex)_n =0$ for all $n> \reg(E)$. 

Each element $a\in A$ has a natural image, denoted by $a^*$, in $G_I(A)$.  Thus, if $x\in I\setminus I^2$, then $x^*$ is a filter-regular  element on $G(\FilM)$ if and only if $(M_{n+2}: x)\cap  M_n = M_{n+1}$ for all  $n> \reg(G(\FilM))$.  
 
 \begin{Lemma} \label{Fil-Property}
 Let $x\in I\setminus \mm I$ be an element such that $x^*$ is a filter-regular on $G(\FilM)$ and let $a= \reg (G(\FilM))$. Then
 \begin{itemize}
 \item[(i)] $r(\FilM) \le a$.
 \item[(ii)] $xM \cap M_n = xM_{n-1}$ for $n \ge a+1$.
 \item[(iii)] $M_{n+1} : x = M_n + (0_M:x)$  and $(0_M:x) \cap M_{n+1} = 0$ for $n\ge a$.
 \end{itemize}
 \end{Lemma}
\begin{pf} As mentioned above, $r(\FilM)$ is the largest generating degree of $G(\FilM)$. Hence (i) is a well-known fact. The statements (ii) and (iii) are just filtration version of \cite[Lemma 4.4]{T2} and \cite[Proposition 4.6]{T2}, respectively. Eventually, (ii) and (iii) characterize the filter-regular property of $x^*$, see \cite[Theorem 1.5]{RV}.
\end{pf}

We are interested in bounding  $\reg(G(\FilM))$.  In the case of $I$-adic filtration, i.e. $M_n = I^nM$ for all $n \geq  0$, this problem was solved by Rossi, Trung and Valla  provided  $M=R$  and $I=\mm$ (see \cite{RTV}) and by Linh in the general case. The bound was given in terms of the so-called  extended degree of $M$ w.r.t. $I$ (see \cite[Theorem 4.4]{L1}). We now recall the definition of this notion.
 
 \begin{Definition} \label{A1}  {\rm (see \cite[Section 3]{L1})
 An {\it extended degree} $D(I,M)$ of $M$ w.r.t. an $\mm$-primary ideal $I$ is a numerical function satisfying the following properties:
 
  (i)  $D(I,M) = D(I,M/L) + \ell (L)$, where $L = H^0_\mm(M)$,

 (ii)  $D(I,M) \geq D(I,M/xM)$ for a generic element $x\in I \setminus  \mm I$ on $M$,

 (iii) $D(I,M) = e(I,M)$ if $M$  is a Cohen-Macaulay $A$-module,
where $e(I ,M)$ denotes the multiplicity of $M$ w.r.t. $I.$ }
 \end{Definition}
 
This notion extends that of extended degrees of graded modules in \cite{DGV} and \cite[Definition 9.4.1]{Va}. Assume that $A$ is a homomorphic image of a Gorenstein ring $S$ with $\dim S = s$. Then, as an example of extended degrees we can take  the homological degree which is defined recursively as follows 
\begin{gather}
\hdeg(I,M) := e(I,M) + \displaystyle \sum_{i=0}^{d-1} {d-1 \choose
i}\hdeg(I,\Ext_S^{s+i+1-d}(M,S)) \label{E:hdeg}
\end{gather}
if $d= \dim M > 0$, and $\hdeg(I,M) =\ell(M)$ if $d= 0$.   This follows from \cite[Theorem 9.4.2]{Va}.

 For filtered modules, Linh's result \cite[Theorem 4.4]{L1} can be modified as follows:

\begin{Theorem}\label{A3}
Let $M$ be a finitely generated $A$-module with $\dim M = d \geq  1$, $\FilM = \{ M_n\}_{n\geq 0}$ a good $I$-filtration of $M$ and $D(I,M)$ an arbitrary extended degree of $M$ with respect to $I$.  Then

\rm{(i)} \ $\reg(G(\FilM) )\leq D(I,M) +r(\FilM) - 1 \ {\mathrm{if}} \  d = 1$,

\rm{(ii)} $\reg(G(\FilM)) \leq [D(I,M) + r (\FilM)+1]^{3(d-1)!-1} - d\  {\mathrm{if}}\  d\geq 2$.
\end{Theorem}

Note that the  bound  in (ii) is  simpler and in most cases  a little bit better than the one  in \cite[Theorem 4.4]{L1}. The proof is similar to that in \cite{L1}, but it needs some modifications. The main reason for the modifications is the fact that unlike the $I$-adic case, the module $G(\FilM)$ is not generated in degree $0$. Therefore below we formulate the required modifications and only give a proof when it is necessary.

 We call
$$H_{\FilM}(n) = \ell(M/M_{n + 1})$$
the Hilbert-Samuel function of $M$ w.r.t  $\FilM$. Since $H_{\FilM}(n) = \ell(M/I^{n + 1 - r}M_r)$ for  all $n \geq  r$, this function agrees with  a polynomial - called the Hilbert-Samuel polynomial  and denoted by $P_{\FilM}(n)$ - for $n \gg 0$. The Hilbert-Samuel  function $\ell(M/I^{n+1}M)$ and its Hilbert-Samuel polynomial of $M$ w.r.t.  an $\mm$-primary ideal $I$ are denoted by $H_{I,M}(n)$ and $P_{I,M}(n)$, respectively.  Whenever $R_0$ is an Artinian local ring, we denote  the Hilbert function $\ell_{R_0}(E_n)$ and the Hilbert polynomial of a finitely generated graded module $E$ over a graded ring $R$ by $h_E(n)$ and $p_E(n)$, respectively.  Note that
\begin{equation} \label{Hh} H_{\FilM}(n) = \sum_{j=0}^n h_{G(\FilM)}(j).\end{equation}

\begin{Lemma}\label{A4} {\rm (cf. \cite[Lemma 3.5]{L1})}
Let $x \in  I \setminus  \mm I$ such that the initial form $x^*$ of $x$ in $G_I(A)$ is a filter-regular element on $G(\FilM)$. Let $N = M/xM$. Then
$$p_{G(\FilM)}(n) \leq  H_{I,N}(n)$$
for $n \geq  \reg(G(\FilM/xM)).$
\end{Lemma}
\begin{pf}
By a filtration version of the so-called Singh's formula (see \cite[Lemma 1.9]{RV}) we have
$$h_{G(\FilM)}(n) = H_{\FilM/xM}(n) - \ell(M_{n+1} : x/M_n).$$
Using Lemma \ref{Fil-Property}(iii) and essentially the same proof of \cite[Lemma 3.5]{L1} we get
$$p_{G(\FilM)}(n) \leq  H_{\FilM/xM}(n)$$
for $n \geq  \reg(G(\FilM/xM))$.  Since $I^{n+1}M \subseteq  M_{n+1}$, we have
$$H_{\FilM/xM}(n)  = \ell(M/M_{n+1} + xM)  \leq \ell(M/I^{n+1}M + xM) = H_{I,N}(n) .$$
\end{pf}

Recall that a subideal $Q\subseteq I$ is said to be a reduction of $I$ if $I^{n+1} = QI^n$ for $n\gg 0$. If $Q$ is a reduction of $I$, which does not properly contain a reduction of $I$, then $Q$ is called a {\it minimal reduction} of $I$. The second statement of the following result is \cite[Theorem 3.6]{L1}. However our proof here is much shorter.

  \begin{Lemma}\label{A5}
 Let $\dim M=d\geq 1$ and $I$ be an $\mm$-primary ideal. Then

{\rm (i)} $\ell(M/I^{n+1}M)\leq \binom{n+d}{d}\ell(M/QM),$ 
where $Q$ is a minimal reduction ideal of $I(A/\Ann(M)),$

{\rm (ii)} $\ell(M/I^{n+1}M)\leq \binom{n+d}{d}D(I,M).$ 
\end{Lemma}
\begin{pf}

(i) Let $Q = (x_1, ..., x_d)$. From the epimorphism 
$$ B := (M/QM)[x_1, ..., x_d] \longrightarrow \bigoplus_{n\geq 0} Q^nM/Q^{n+1}M,$$
we get that
$$\ell(M/I^{n+1}M) \leq  \ell(M/Q^{n+1}M) \leq  \sum_{i=0}^{i=n} \ell(B_i) \leq \binom{n+d}{d}\ell(M/QM).$$

(ii) We may choose $x_1, ..., x_d \in  I$ such that $x_i$ is a generic element on 
$M/(x_1, ..., x_{i-1})M$. Then $Q = (x_1, ..., x_d)$ 
is a minimal reduction of $I(A/\Ann(M))$. 
By (ii) and (iii) of Definition \ref{A1}, we get
$$D(I,M) \geq  D(I,M/(x_1, ..., x_d)M) = \ell(M/(x_1, ..., x_d)M) = \ell(M/QM).$$
Hence the statement follows from (i).\end{pf}

Recall that the geometric Castelnuovo-Mumford regularity of a graded $R$-module $E$ is
$$\geom(E) = \max \{a_i(E) + i \mid  i \geq  1 \}$$
(see \cite{RTV}). The following result is a module version of \cite[Corollary 5.3]{H} and the proof is essentially the same   (see also \cite[Proposition 3.5]{HZ} and \cite[Lemma 4.2]{L1}).
\begin{Lemma} \label{A6} 
Let $M$ be a finitely generated $A$-module such that $\depth M > 0.$ Then
$$\reg(G(\FilM)) = \geom(G(\FilM)).$$
\end{Lemma}

The next result is an extension of  \cite[Lemma 4.3]{L1} to the filtration case.

\begin{Lemma}\label{A7} 
Let $\overline{M} = M/H^0_\mm(M)$. Denote the filtration $\FilM/ H^0_\mm(M)$ of $\overline{M}$ by $\overline{\FilM}$. Then
$$\reg(G(\FilM)) \leq  \max \{ \reg(G(\overline{\FilM}));\ r(\FilM)\} + \ell(H^0_\mm(M)).$$
\end{Lemma}
\begin{pf}
Let $L = H^0_\mm(M)$ and 
$$K = \bigoplus_{n\geq 0}(M_{n+1} + M_n \cap  L)/M_{n+1} \cong  \bigoplus_{n \geq  0} \frac{M_n \cap L}{M_{n+1}\cap L}.$$
By the Artin-Rees theorem there is an integer $c$ such that 
$$M_{n+1} \cap  L \subseteq  I^{n+1-r}M \cap  L \subseteq  I^{n+1-r-c}L = 0$$
for $n \gg 0.$ Hence, $\ell(K) = \ell(L)$. It is clear that we have an exact sequence  
\begin{gather}  
0 \longrightarrow  K \longrightarrow  G(\FilM) \longrightarrow  G(\overline{\FilM}) \longrightarrow  0.
\label{E:A7}
 \end{gather}
 Let $p = \max \{ \reg(G(\overline{\FilM}));\  r \}$. Then there is an integer $p + 1 \leq m \leq p + \ell(L) + 1$ such that $K_m = 0$.  Since $m >  \reg(G(\overline{\FilM}))$, from (\ref{E:A7}) it then implies that $H^i_{G_I(A)_+}(G(\FilM))_{m-i} = 0$ for all $i \geq  0$. Note that $G(\FilM)$ is generated over $G_I(A)$ by elements of degrees $\leq  r \leq  m -1$. Hence, by \cite[Lemma 2.1]{Na}   $\reg(G(\FilM)) \leq  m -1 \leq  p + \ell(L)$, as required.  
\end{pf}

We can  now prove  Theorem \ref{A3}. As  said above, it is essentially the same as the proof of \cite[Theorem 4.4]{L1}. Hence, in the proof  below we only give  the details when  modifications are needed.

\vskip0.5cm
\noindent {\it Proof of Theorem \ref{A3}}.
Let $G=G_I(A)$, $L = H^0_\mm(M)$ and $x \in  I \setminus  \mm I$ be a generic element on $M$. Let $\overline{M} = M /  L$ and recall that $ r : = r(\FilM)$. By Lemma \ref{A7},
$$\reg(G(\FilM)) \leq  \max \{ \reg(G(\overline{\FilM}));\ r\} + \ell(L).$$
By (i) of Definition \ref{A1}, $D(I,M) = D(I,\overline{M}) + \ell(L)$. Hence, we only need to show the following statements:

(i')  $\reg(G(\overline{\FilM})) \leq D(I,\overline{M})  + r - 1 \ {\mathrm{if}} \  d = 1$,

(ii') $\reg(G(\overline{\FilM})) \leq [D(I,\overline{M}) + r + 1]^{3(d-1)!-1} - d\  {\mathrm{if}}\  d\geq 2$.

\noindent Replacing $M$ by $\overline{M}$, we may assume that $\depth M > 0$.  Set $D= D(I,M)$.  By Lemma \ref{A6},
\begin{equation}
\reg(G(\FilM)) = \geom (G(\FilM)).\label{E:A3}
\end{equation}

If $d = 1$ then $M$ is a Cohen-Macaulay module, $G(\FilM)$ is a $G$-module of dimension one generated by elements of degree at most $r$. Hence, by Lemma \ref{A6} and  \cite[Lemma 2.2]{L1}, we have
$$\reg(G(\FilM)) = \geom(G(\FilM)) =a_1(G(\FilM)) + 1 \leq  e(G(\FilM))  + r - 1 = e(I,M) +r - 1 .$$
The last equality follows from \cite[Proposition 11.4(iii)]{AM}. Hence, by Definition \ref{A1}(ii),  $\reg(G(\FilM)) \le D + r-1$. 

Note that $r(\FilM/ xM) \le r$ and $D >0$. If $d \geq  2$, let $N = M/xM$ and $m: = \max \{ r;\  \reg(G(\FilM/xM)) \}$. Using the following exact sequence
\begin{equation} \label{Filsq} 0 \rightarrow  \bigoplus_{n\ge 1}\frac{xM \cap M_n}{xM_{n-1} + xM\cap M_{n+1}} \rightarrow G(\FilM)/x^*G(\FilM) \rightarrow G(\FilM/xM) \rightarrow 0,
\end{equation}
 and Lemma \ref{Fil-Property}(ii) we get $\geom(G(\FilM)/x^*G(\FilM)) =  \geom(G(\FilM/xM))$. Hence,
 $$\geom(G(\FilM)/x^*G(\FilM)) \leq  m.$$
Since  $G(\FilM)$ is  generated by elements of degree at most $r \leq  m$, by \cite[Theorem 2.7]{L1},
$$\geom(G(\FilM)) \leq  m + p_{G(\FilM)}(m).$$
By Lemma \ref{A4} and Lemma \ref{A5} (ii) and the fact that $D(I,N) \leq  D$, we get
\begin{equation}
\geom(G(\FilM)) \leq  m + D(I,N){m+d-1\choose d-1} \leq m+ D{m+d-1\choose d-1}.\label{E:A3b}
\end{equation}
If $d = 2$, then by (i) of the theorem 
$$m = \max \{ r;\  \reg(G(\FilM/xM)) \} \leq  D(I,N) + r - 1 \leq D + r  -1.$$
Since $r\ge 0$, by (\ref{E:A3}) and (\ref{E:A3b}), we get
\begin{align*}
\reg(G(\FilM)) = \geom(G(\FilM)) & \leq  m + D(m+1) \\
                                                    & \leq  D^2 +(r + 1)D + r - 1 \\
                                                    & \leq  (D + r + 1)^2 - 2.
\end{align*}
Let $d \geq  3$. The case $m=0$ is trivial, so we may assume $m>0$. From (\ref{E:A3b}) we get
\begin{gather}\geom(G(\FilM)) \leq  m + D{m+d-1 \choose d-1} \leq  D(m + 1)^{d - 1} -1. \label{E:A3c}
\end{gather}
By the induction hypothesis we may assume that 
$$ m \leq  [D(I, N) + r  + 1]^{3(d-2)! -1} -  d+1\leq (D + r + 1)^{3(d-2)! - 1} - d+ 1.$$
Hence, by (\ref{E:A3}) and (\ref{E:A3c}),
$$\reg(G(\FilM)) \leq  [D + r + 1]^{3(d-1)! - 1} - d.$$
\ \hfill $\square$

\vskip0.5cm
If we write
$$P_\FilM (t) = \sum_{i=0}^d (-1)^i e_i(\FilM){t+d-i \choose d-i},$$
then the integers $e_i(\FilM)$ are called {\it Hilbert coefficients} of $\FilM$ (see \cite[Section 1]{RV}). Note that $e_0(\FilM),...,e_{d-1}(\FilM)$ are the Hilbert coefficients of the graded module $G(\FilM)$. In the last section we need an estimation of  Hilbert coefficients. This was done for the $\mm$-adic filtration of a ring in \cite[Theorem 4.1]{RTV} and extended to the module case in \cite[Theorem 3.1]{L2}. However the proof of \cite[Theorem 3.1]{L2} has a gap. Namely,  in the case $d=1$ there was used the following wrong inequality: $|(e+1)e(I,M) - \ell(M/I^{r+1}M)| \le |(e+1)e(I,M) - (r+1)|$. Therefore we give here the proof of the following result, which moreover improves \cite[Theorem 3.1]{L2}. 

\begin{Theorem} \label{Hilb} Let $M$ be a finitely generated $A$-module with $\dim M = d \geq  1$ and $\FilM = \{ M_n\}_{n\geq 0}$ a good $I$-filtration of $M$ . Then

\item[(i)] $e_0(\FilM) = e(I,M) \le D(I,M)$;

\item[(ii)] $|e_1(\FilM) | \le (D(I,M) + r(\FilM) -1)D(I,M)$;

\item[(iii)] $| e_i(\FilM)|\le (D(I,M) + r (\FilM) +1)^{3i! -i +1} $ if $i\ge 2$. 
\end{Theorem}

\begin{pf}   (i): By \cite[Proposition 11.4(iii)]{AM} $e_0(\FilM) = e(I,M)$. By Definition \ref{A1}, $e(I,M) \le D(I,M)$.

(ii) -(iii): From the Grothendieck-Serre formula
$$h_{G(\FilM)}(n) - p_{G(\FilM)}(n) = \sum_{j=0}^d(-1)^j H^j_{G_+}({G(\FilM)})_n,$$
it follows that $h_{G(\FilM)}(n) = p_{G(\FilM)}(n)$ for all $n>\reg(G(\FilM))$.
By (\ref{Hh}) it follows that
\begin{equation} \label{Ph} \ell(M/M_{m+1}) = \sum_{i=0}^d (-1)^i e_i(\FilM){m+d-i \choose d-i}
\end{equation}
for any $m\ge \reg(G(\FilM))$. For  simplicity we set $r= r(\FilM)$, $D:= D(I,M)$ and $e_i := e_i(\FilM)$.

Assume that $d=1$. Using  Theorem \ref{A3}(i) and putting $m = D+r -1$ into the equality (\ref{Ph}), we have
$$e_1 =(D+r) e_0 - \ell(M/M_{D+r}).$$
Since $M_{n} = I^{n-r}M_r$ for $n\ge r$ and $M_r \neq 0$, 
$$\ell(M/M_{D+r})\ge \ell(M_r/IM_r) + \cdots +  \ell(I^{D-1}M_r/I^D M_r) \ge D.$$
By the first claim, this implies 
$$e_1 \le  (D+r)e_0 - D \le D (D+r) - D= D(D+r-1).$$
On the other hand, since $r\ge 0$, by  Lemma \ref{A5}(ii), 
$$- e_1 =- (D+r) e_0 +  \ell(M/M_{D+r})\le D(D+r) - (D+r) \le D(D+r-1).$$
Hence $|e_1| \le  D(D+r-1)$ and the case $d=1$ is proven. 

Let $d\ge 2$. First we assume that $\depth M >0$. Let $x\in I\setminus \mm I$ be a generic element. Then $x^*\in G$ is a filter-regular element on $G(\FilM)$.  By
\cite[Proposition 1.10]{RV}, $e_i = e_i(\FilM/xM)$ for all $i<d$. Note that $0 \le r(\FilM/xM) \le r$ and by Definition \ref{A1}(ii),  $D(I, M/xM) \le D$. Using  the induction hypothesis we then get
\begin{equation}\label{EHilb0}  |e_1| \le D(D+r- 1) \ \ \text{and} \ \ | e_i| \le (D + r +1)^{3i! -i +1} \ \text{for} \ 2\le i \le d-1.\end{equation}
To prove the inequality for $e_d$, we set $\mu = (D+r+1)^{3(d-1)!-1}$. By Theorem \ref{A3}, $\reg(G(\FilM) ) \le \mu -d$. Since $\reg(G(\FilM)) \ge r\ge 0$, $\mu \ge d$. Hence putting $m= \mu -d$ into  (\ref{Ph}), we have
\begin{eqnarray}\nonumber  |e_d| & = |(\ell(M/M_{\mu -d +1}) - e_0{\mu -d  +d \choose d}) + \sum_{i=1}^d (-1)^i e_i{\mu -d +d-i \choose d-i}| \\
&\nonumber  \le  |\ell(M/M_{\mu -d +1}) - e_0{\mu    \choose d}| + \sum_{i=1}^{d-1} | e_i|{\mu  -i \choose d-i} \\
& \le \max\{\ell(M/M_{\mu -d +1}),  \ e_0{\mu    \choose d}\} + \sum_{i=1}^{d-1} | e_i|{\mu  -i \choose d-i}. \label{EHilb1}
\end{eqnarray}
Note that  ${\mu    \choose d} \le \mu^d$. By (i) and  Lemma \ref{A5}(ii) it yields
\begin{equation} \label{EHilb2}  \max\{\ell(M/M_{\mu -d +1}),  \ e_0{\mu    \choose d}\} \le D\mu^d. \end{equation}
Further, by (\ref{EHilb0})
\begin{equation} \label{EHilb3} 
|e_1|{\mu  -1 \choose d-1} \le D(D+r-1)\mu^{d-1},
\end{equation}
 and
 \begin{equation} \label{EHilb4} 
 \sum_{i=2}^{d-1} | e_i|{\mu  -i \choose d-i} \le  \sum_{i=2}^{d-1} (D+r+1)^{3i! -i +1}\mu^{d-i} \le \mu^{d-1}\sum_{i=0}^{d-2}\frac{1}{2^i} < 2\mu^{d-1}.
\end{equation}
Since $D(D+r- 1) + 2<  (D+r+1)^2 \le \mu $, from (\ref{EHilb1}) - (\ref{EHilb4}) we obtain
$$|e_d| \le D\mu^d + \mu^d =  (D+1)(D+r+1)^{3d! -d} \le (D+r+1)^{3d! - d+1}.$$
Finally assume that $\depth M =0$. Set $L = H^0_\mm(M)$, $\bar{M} = M/L$  and $\overline{\FilM} = \FilM/L$. Then
$$\begin{array}{ll} \ell(\frac{M}{M_{n+1}}) & = \ell(\frac{M}{M_{n+1}+L}) + \ell(\frac{M_{n+1} + L}{M_{n+1}}) =  \ell(\frac{M}{M_{n+1}+L}) + \ell(\frac{L}{L\cap M_{n+1}})\\
& =  \ell(\frac{M}{M_{n+1}+L})  +\ell(L)
\end{array}$$
for $n\gg 0$ (by Lemma \ref{Fil-Property}(iii)). Hence $e_i = e_i(\overline{\FilM})$ for $i\le d-1$ and $e_d = e_d(\overline{\FilM}) +(-1)^d\ell(L)$. Since $D = D(I, \bar{M}) + \ell(L)$ and $r(\overline{\FilM}) \le r$, we get
$$\begin{array}{ll}
 |e_1|  & = e_1(\overline{\FilM}) \le D(I,\bar{M}) (D(I,\bar{M})+r-1) \le  D(D+r-1)\\
 | e_i|  & = e_i(\overline{\FilM}) \le (D(I, \bar{M}) + r +1)^{3i! -i +1} \le (D + r +1)^{3i! -i +1} \ \text{for} \ 2\le i \le d-1,\\
 | e_d|  & \le e_d(\overline{\FilM})  +\ell(L) \le (D(I, \bar{M}) + r +1)^{3d! -id+1}  + \ell(L) \le (D + r +1)^{3d! -d+1}.
 \end{array}$$
\end{pf}

\section{Regularity of associated graded rings: Graded case} \label{Graded}
Let $A = \oplus _{n \geq  0}A_n$ be a Noetherian standard  graded algebra over an Artinian local ring $(A_0, \mm_0)$, i.e. $A = A_0[A_1]$. As usual, we assume that $A_0 / \mm_0$ is infinite.  We denote the maximal homogeneous  ideal $\mm_0 \oplus (\oplus_{n \geq  1} A_n)$ of $A$  by $\mm$. Let   $M = \oplus_{n\in \ZZ}M_n$  be a finitely generated graded $A$-module of dimension $d$ and $\FilM = \{ \Mc_n\}$  a good $I$-filtration consisting of homogeneous submodules of $M$, where $I$ is a homogeneous $\mm$-primary  ideal of $A$. Our goal is to bound $\reg(G(\FilM))$ in terms of $\reg(M)$ and  $r(\FilM)$ in the case $A$ is a polynomial ring over a field.  

 More generally, let $A$ be a homomorphic image of a Gorenstein graded algebra $S$ with $\dim S =s$. We define $e(I,M)$  to be the degree of the graded module $G_I(M)$, or equivalently, $(d-1)!$ times of the leading coefficient of the Hilbert polynomial $h_{G_I(M)}(t)$. Then, as in the local case, one can define  the homological degree $\hdeg(I,M)$  by  formula (\ref{E:hdeg}). In particular, if $I = \mm$, we set 
\begin{gather}
\hdeg(M) := \hdeg(\mm,M)= e(M) + \displaystyle \sum_{i=0}^{d-1} {d-1 \choose
i}\hdeg(\mm, \Ext_R^{s  + i +1 -d}(M,R)).\label{E:hdeg2}
\end{gather}
In fact, this definition was first given in \cite{DGV} and \cite[Definition 9.4.1]{Va}. 

We now bound $\reg(G(\FilM))$ in terms of  $r(\FilM)$ and $\hdeg(M)$. Our idea is to reduce to the local case.

\begin{Lemma}\label{B1}
Let $I \subset  A$ be a homogeneous $\mm$-primary ideal. Then
$$\hdeg(I_\mm,M_\mm) \leq  \ell(A/I)^d\hdeg(M).$$
\end{Lemma}

\begin{pf}
Let $p = \ell(A/I) = \ell((A/I)_\mm)$. There is a composition series
$$0 = L_0 \subset L_1 \subset  ... \subset  L_p = (A/I)_\mm.$$
Hence, $\mm^p A_\mm \subseteq I_\mm$. By \cite[Lemma 1.3]{L2}, we have 
$$\hdeg(I_\mm,M_\mm) \leq  p^d \hdeg(\mm A_\mm,M_\mm).$$
Since the $\Ext$ functor  commutes with the localization (see, e.g. \cite[Theorem 9.50]{Ro}), from the recursive  formulas (\ref{E:hdeg}) and (\ref{E:hdeg2}) and the fact that $e(E) = e(E_\mm)$ for any finitely generated graded $A$-module $E$, it follows that $\hdeg(\mm A_\mm,M_\mm) = \hdeg(M)$. Hence, $\hdeg(I_\mm,M_\mm) \leq  p^d \hdeg(M).$
\end{pf}

\begin{Theorem}\label{B3}
Let $\FilM$ be a good $I$-filtration of  a graded $A$-module $M$  of dimension $d$.  Then

\rm{(i)}  $\reg(G(\FilM)) \leq \ell(A/I)\hdeg(M) + r(\FilM) - 1 \ {\mathrm{if}} \  d = 1$,

\rm{(ii)} $\reg(G(\FilM)) \leq [\ell(A/I)^d\hdeg(M) + r (\FilM) +1]^{3(d-1)!-1} -d\  {\mathrm{if}}\  d\geq 2$.
\end{Theorem}
\begin{pf} Denote by $\RR = A \oplus  I \oplus I^2  \oplus  ...$ the Rees algebra of $A$ w.r.t. $I$ and $\RR_+ = \oplus _{n\geq 1} I^n$. Then we can consider $G(\FilM)$ as a finitely generated module over $\RR$. If $E = \oplus_{n \in  \ZZ} E_n$ is a graded module over $\RR$, we denote by $E_\mm$ and $(E_n)_\mm$ the localization of $E$ and $E_n$ as $A$-modules with respect to the multiplicative set $A \setminus  \mm$.  We can consider $G_I(A)$ and $G(\FilM)$ as graded modules over $\RR$. Then it is easy to see that $(G(\FilM))_\mm \cong  G(\FilM_\mm)$ and $(G_I(A))_\mm \cong G_{I_\mm}(A_\mm)$. This implies
 $$\reg(G(\FilM)) = \reg(G(\FilM_\mm)).$$
On the other hand, by Lemma \ref{B1},
$$\hdeg(I_\mm,M_\mm) \leq  \ell(A/I)^d\hdeg(M).$$
We also have $r(\FilM_\mm) \leq  r(\FilM).$ Hence, applying Theorem \ref{A3} to the homological degree we get the statements of the theorem.
\end{pf}

An important consequence of Theorem \ref{B3} is
\begin{Corollary}\label{B4}
Let $I$ be an $\mm$-primary homogeneous ideal of a polynomial ring $A = K[x_1,...,x_n]$ over an infinite field $K$. Then
 
\rm{(i)}  $\reg(G_I(A)) \leq \ell(A/I) - 1 \ {\mathrm{if}} \  d = 1$,

\rm{(ii)} $\reg(G_I(A)) \leq (\ell(A/I) +1)^{3(d-1)!-1} - d\  {\mathrm{if}}\  d\geq 2$.
\end{Corollary}

If  $M$ is an arbitrary graded module over a polynomial ring $A$, then we can bound $\reg(G(\FilM))$ in terms of $\reg(M),\ r(\FilM)$ and some other invariants of $M$ as follows:

\begin{Theorem}\label{B5}
Let $M$ be a finitely generated graded module of dimension $d$ over a polynomial  ring $A = K[x_1,...,x_n]$. Let $i(M)$ denote the initial  degree of $M$ (i.e. $i(M) = \min\{ p\mid  M_p \ne 0 \}$) and $\mu(M)$ the minimal number of generators of $M$. Then

\rm{(i)}  $\reg(G(\FilM)) \leq \ell(A/I)\mu (M)[\reg(M) - i(M) +1]^n +r(\FilM) - 1 \ {\mathrm{if}} \  d = 1$,

\rm{(ii)} $\reg(G(\FilM)) \leq [\ell(A/I)^d (\mu (M)(\reg(M) - i(M) +1)^n)^{2^{(d - 1)^2}} + r(\FilM)+1]^{3(d - 1)! - 1} - d\  {\mathrm{if}}\  d\geq 2$.
\end{Theorem}
\begin{pf}
By \cite[Theorem 5.1]{CHH},
\begin{align*}
\hdeg(M) &\leq  \biggl[\mu(M) \binom{\reg(M) - i(M) +n}{n}\biggr]^{2^{(d-1)^2}} \\
                & \leq [\mu(M) (\reg(M) - i(M) +1)^n]^{2^{(d - 1)^2}}.
\end{align*}
Hence, the statement follows from Theorem \ref{B3}.
\end{pf}

In the rest of  this section  we assume that $\FilM$ is a good $I$-filtration, where $I$ is an $\mm$-primary ideal  generated by homogeneous elements of the same degree. Under this setting we will give a bound on $\reg(G(\FilM))$ which is much better than the ones in Theorem \ref{B5} and holds for any standard graded ring $A$. In this new setting, we have
\begin{Lemma}\label{C1}
Let $i(M)$ denote the initial degree of $M$. Let $\overline{M} := M/H^0_\mm(M)$ and $\overline{\FilM} :=  \FilM/H^0_\mm(M)$. Then
$$ \reg(G(\FilM)) \leq  \max\{\reg(G(\overline{\FilM})); \ \reg(M) -i(M) +r(\FilM) \}.$$
\end{Lemma}
\begin{pf}
As in the proof of the Lemma \ref{A7}, we have the following exact sequence
\begin{gather}
0\longrightarrow K\longrightarrow  G(\FilM) \longrightarrow  G(\overline{\FilM})  \longrightarrow  0,\label{E:C1}
\end{gather}
where
$$K = \bigoplus_{n \geq  0} \frac{\Mc_{n+1}+\Mc_n\cap  H^0_\mm(M)}{\Mc_{n+1}}$$
is a module of finite length. 

Let $n \geq  \reg(M) - i(M) + r + 1$, where  $r := r(\FilM)$.  Note that $\reg(M) \geq  i(M)$ and $I$ is generated by homogeneous elements of degree at least one. Hence, 
$$ \Mc_n = I^{n-r}\Mc_r \subseteq I^{n-r}M \subseteq  \bigoplus _{p \geq  \reg(M) + 1}M_p.$$
 Since $H^0_\mm(M) \subseteq  M$ and $H^0_\mm(M)_p = 0$ for all $p \geq  \reg(M) + 1$, we get that $\Mc_n \cap H^0_\mm(M) = 0$. Hence, $K_n = 0$. From this fact  and the exact sequence (\ref{E:C1}) we get the statement of the lemma.
\end{pf}

\begin{Lemma}\label{C2}
Let $x \in  I \setminus  \mm I$ be a homogeneous element. Assume that the initial form $x^*$ of $G_I(A)$ is a filter-regular element on $G(\FilM)$. Then $x$ is a filter-regular element on $M$.
 \end{Lemma}
 
\begin{pf}
The element $x^*$ is a filter-regular element on $G(\FilM)$ means that
\begin{gather}
(\Mc_{n+2}:x) \cap  \Mc_n = \Mc_{n + 1}\label{E:C2}
\end{gather}
for all $n \geq  n_0$, where $n_0$ is a certain fixed number. Let $u \in  (0:_Mx)_p$. Since $M$ is finitely generated and $M/\Mc_{n_0}$ is of finite length, it follows that $u \in  \Mc_{n_0}$ when $p \gg 0$. Then, by (\ref{E:C2}),
$$ u \in  (0:_Mx) \cap  \Mc_{n_0} \subseteq  (\Mc_{n_0+2}:x) \cap  \Mc_{n_0} = \Mc_{n_0+1}.$$
By induction it implies that
$$ u \in  \bigcap_{q \geq  n_0} \Mc_q =  \bigcap_{q \gg 0} I^{q-r}\Mc_r = 0.$$
This means $(0:_Mx)_p = 0$ for $p \gg 0$, or equivalently $x$ is a filter-regular element on $M$.
\end{pf}

We can now improve Theorem \ref{B5} in the case of an equi-generated ideal $I$ as follows:
\begin{Theorem}\label{C3}
Assume that $I$ is generated by elements  of degree $\Delta \ge 1$. Let $Q$ be a homogeneous minimal reduction of $I(A/\Ann(M))$. Let  $i(M)$ denote the  initial degree of $M$. Then

\rm{(i)} \ $\reg(G(\FilM)) \leq \ell(M/QM) + r (\FilM) + \reg(M) - i(M)    - 1 \ {\mathrm{if}} \  d = 1$,

\rm{(ii)} $\reg(G(\FilM)) \leq [\ell(M/QM) + r(\FilM) + \reg(M) -i(M) + (d-1)\Delta ]^{3(d-1)!-1} - d\  {\mathrm{if}}\  d\geq 2$.
\end{Theorem}

\begin{pf} The main idea is the same as in the proof of \cite[Theorem 4.4]{L1}.
By Lemma \ref{C1}, it suffices to consider the case $\depth M > 0$. Let $r:= r(\FilM)$.

(i) If $d = 1$, then $M$ is a Cohen-Macaulay module. Hence, by \cite[Lemma 2.2]{L1},  Lemma \ref{A6} and   \cite[Proposition 11.4(iii)]{AM}, we get
\begin{align*}
\reg(G(\FilM)) & = \geom(G(\FilM)) = a_1( G(\FilM))  + 1\leq  e(G(\FilM))  + r -1 \\
                       & =  e(I,M) + r-1 = e(Q,M) +r - 1 \leq \ell(M/QM) +r  - 1.
\end{align*}

(ii) $d \geq  2$. It is clear that $Q$ is generated by $d$ elements of degree $\Delta $. Then one can find a minimal basis $\{x_1, ..., x_d \}$ of $Q$ such that the initial form $x^*_1$ in $G$ is a filter-regular element on $G(\FilM)$ (see \cite[Lemma 3.1]{T}). Note that all elements $x_1,...,x_d$ have degree $\Delta  \geq  1$. Let $x = x_1$. Then $N:= M/xM$ is again a graded module. Let $m = \max\{r; \reg(G(\FilM/xM)) \}$. Using the exact sequence (\ref{Filsq}) we then get $\geom(G(\FilM)/x^*G(\FilM)) = \geom(G(\FilM/xM)) \le m$. Note that $(x_2, ... x_d)$ is a minimal reduction of $I(A/\Ann(N))$. Hence, by  \cite[Theorem 2.7]{L1},  Lemma \ref{A4} and Lemma \ref{A5} (i) we get
\begin{align}
\geom(G(\FilM))  &\leq m + \ell(N/(x_2, ..., x_d)N) \binom{m + d - 1}{d-1}\notag \\
                           &\leq m+ \ell(M/QM)\binom{m + d - 1}{d-1},\label{E:C4}
\end{align}
which is similar to (\ref{E:A3b}).
Note that $\reg(N) \leq  \reg(M) + \Delta  -1, \ r (\FilM/xM)\le r$ and $i(N) \geq i(M)$. Let $d = 2$. By the induction hypothesis we have
\begin{align*}
\reg(G(\FilM/xM)) & \leq \ell(N/(x_2,...,x_d)N) + r + \reg(N) -i(N) -1 \\
                    & \leq  \ell(M/QM) + r+\reg(N)-i(M)-1\\
                    & \le  \ell(M/QM) +r+\reg(M)-i(M)+\Delta  - 2.
\end{align*}
Hence, 
$$m \le \ell(M/QM) +r+\reg(M)-i(M)+\Delta -2 .$$
Together  with (\ref{E:C4}) and Lemma \ref{A6} this yields 
$$\reg(G(\FilM)) = \geom(G(\FilM)) \leq [\ell(M/QM)+r+\Delta  +\reg(M) -i(M)]^2-2.$$
If $d\geq 3$, then again by the induction  hypothesis we get
\begin{align*}
m & \leq  [\ell(N/(x_2, ... x_d)N)+r+(d-2)\triangle+\reg(N)-i(N)]^{3(d-2)!-1} -d+1\\
     & \leq  [\ell(M/QM)+r+(d-1)\triangle+\reg(M)-i(M)]^{3(d-2)!-1} -d+1.
\end{align*}
Hence, using (\ref{E:C4}) we can complete the proof of (ii).
\end{pf}

\section{Regularity of fiber cones} \label{Fiber}

In this section we assume that $\FilM = \{ M_n\}$ is a good $I$-filtration of a finitely generated module $M$ over a local ring $(A,\mm)$, where $I$ is an $\mm$-primary ideal. Given an ideal $\qq$  containing $I$, we define
$$F_\qq(\FilM) := \oplus_{n\ge 0}M_n/\qq M_n,$$
and call it the {\it fiber cone} of $\FilM$ with respect to $\qq$. This notion was introduced in \cite[Section 5]{RV} (see also  \cite{JV}). If $\FilM$ is the $I$-adic filtration of $A$ and $\qq=\mm$, then this is  the classical fiber cone $F_\mm(I) = \oplus_{n\ge 0}I^n/\mm I^n$ of $I$. Note that $F_\qq(\FilM)$ is a graded module over $G = G_I(A)$. 

The purpose of this section is to  give a bound for the Castelnuovo-Mumford regularity $\reg (F_\qq(\FilM))$ in terms of $D(I,A)$ and $r(\FilM)$. We will apply the results of Section \ref{Local}. Following \cite[(3)]{RV} we define a new good $I$-filtration
$$\qq M:\ M \supseteq \qq M \supseteq \qq M_1 \supseteq \cdots \supseteq \qq M_n \supseteq \cdots $$

\begin{Lemma}\label{F1}
Let $M$ be a finitely generated $A$-module with $\dim M = d \geq  1$, $\FilM = \{ M_n\}_{n\geq 0}$ a good $I$-filtration of $M$ and $D(I,M)$ an arbitrary extended degree of $M$ with respect to $I$.   Assume that $I\subseteq \qq$ and $M_{n+1} \subseteq \qq M_n$ for all $n\ge 0$. Then

\rm{(i)} \ $a_0(F_\qq(\FilM) )\leq D(I,M) +r (\FilM) \ {\mathrm{if}} \  d = 1$,

\rm{(ii)} $a_0(F_\qq(\FilM)) \leq (D(I,M) + r(\FilM) +2)^{3(d-1)!-1} - d\  {\mathrm{if}}\  d\geq 2$.
 \end{Lemma}
 \begin{pf} Since $M_{n+1} \subseteq \qq M_n$ for all $n\ge 0$, by \cite[Proposition 5.1]{RV} we have the following exact sequence of $G$-graded modules
 $$ 0 \rightarrow F_\qq(\FilM) \rightarrow G(\qq \FilM) \rightarrow N(-1) \rightarrow 0,$$
where $N = \oplus_{n\geq 0}\qq M_n/M_{n+1}$. Therefore 
$$ a_0(F_\qq(\FilM) )\leq a_0(G(\qq \FilM)) \leq \reg(G(\qq\FilM)).$$
Note that $r(\qq \FilM) \le r(\FilM) +1 $.  Hence the claim now follows from Theorem \ref{A3}.
\end{pf}

The study of Hilbert coefficients of $\FilM$ in Section \ref{Local} allows us  to bound the Hilbert coefficients of the fiber cone $F_\qq(\FilM)$. 

\begin{Proposition}\label{F2}
Under the assumption of Lemma \ref{F1} we have 
\begin{itemize} 
\item[(i)] $e_0(F_\qq(\FilM)) \leq 2D(I,M)(D(I,M) + r(\FilM) ).$
\item[(ii)] $|e_i (F_\qq(\FilM)) | \leq 2 (D(I,M) + r(\FilM) +2)^{3(i+1)! -i }$ if $1\le i \le d-1$.
\end{itemize}
 \end{Proposition}
 \begin{pf}
 It was shown in \cite[(24)]{RV} that
 $$e_i(F_\qq(\FilM)) = e_i(\FilM) + e_{i+1}(\FilM) - e_{i+1}(\qq \FilM),$$
 for all $0\le i \le d-1$. Let $r:= r(\FilM)$. Since $r(\qq \FilM) \le r+1$, by Theorem \ref{Hilb}, we get
 $$\begin{array}{ll}
 e_0(F_\qq(\FilM))  &\le |e_0(\FilM)| + |e_1(\FilM) | + |e_1(\qq \FilM)|\\
 &\le D+ D(D + r-1) + D(D+r) =  2D(D + r ),
 \end{array}$$
 where $D= D(I,M)$, and for $1\le i\le d-1$:
 $$\begin{array}{ll}
| e_i(F_\qq(\FilM)) | &\le |e_i(\FilM)| + |e_{i+1}(\FilM) | + |e_{i+1}(\qq \FilM)|\\
 &\le (D + r+1)^{3i! - i +1} + (D+r+1)^{3(i+1)! - i} + (D + r+2)^{3(i+1)! - i} \\
 &\le  2(D + r+2)^{3(i+1)! - i}.
 \end{array}$$
 \end{pf}

 \begin{Theorem}\label{F3}
Let $M$ be a finitely generated $A$-module with $\dim M = d \geq  1$, $\FilM = \{ M_n\}_{n\geq 0}$ a good $I$-filtration of $M$.  Assume that $I\subseteq \qq$ and $M_{n+1} \subseteq \qq M_n$ for all $n\ge 0$. Then
\begin{itemize}
\item[(i)] $\reg(F_\qq(\FilM)) \le 2D(I,M)(D(I,M) + r(\FilM) ) + r(\FilM) -1$ if $d=1$;
\item[(ii)] $\reg(F_\qq(\FilM)) \le (D(I,M) +r(\FilM)+2)^2 + D(I,M)^2 - 3$ if $d= 2$;
\item[(iii)] $\reg(F_\qq(\FilM)) \le (D(I,M) +r(\FilM)+2)^{3(d-1)! - 1} - d$ if $d\ge 3$.
\end{itemize}
 \end{Theorem}
 \begin{pf} We do induction on $d$.  Set $D= D(I,M)$ and $r=r(\FilM)$. Let $d=1$. By \cite[Lemma 2.2]{L1} and Proposition \ref{F2}(i) we have
 $$ a_1(F_\qq(\FilM)) + 1 \le e_0(F_\qq(\FilM)) + r-1 \le  2D(D+r ) + r-1.$$
 By Lemma \ref{F1} it yields
 $$\begin{array}{ll} \reg(F_\qq(\FilM)) & = \max\{ a_0(F_\qq(\FilM));\ a_1(F_\qq(\FilM)) + 1\} \\
 &\le \max\{ D+r; \ 2D(D + r )+r-1 \} = 2D(D+r ) + r-1.
 \end{array}$$
 
 Now let $d\ge 2$. Note that both $G(\FilM)$ and $F_\qq(\FilM)$ are modules over $G = G_I(A)$. Hence there exists a generic element $x\in I\setminus \mm I$ such that $x^*\in G$ is a filter-regular element  on $G(\FilM)$ as well as a  filter-regular element  on $F_\qq(\FilM)$ (cf.  \cite[Proposition 2.2]{RV}). Then
 $$F_\qq(\FilM) / x^* F_\qq(\FilM) \cong  \frac{M}{\qq M} \oplus (\oplus_{n\ge 0} \frac{M_n}{\qq M_n + xM_{n-1}}) ,$$
and $$ F_\qq(\FilM/xM) = \oplus_{n\ge 0} \frac{M_n}{\qq M_n + xM \cap M_n}.$$
 Hence we have an exact sequence of $G$-modules:
 $$0 \rightarrow K \rightarrow F_\qq(\FilM)/ x^*F_\qq(\FilM) \rightarrow F_\qq(\FilM/xM) \rightarrow 0,$$
 where
 $$K = \oplus_{n\ge 1} \frac{\qq M_n + xM \cap M_n}{\qq M_n + x M_{n-1}}.$$
 By Lemma \ref{Fil-Property}(ii) $xM \cap M_n = xM_{n-1}$ for all $n> \reg (G(\FilM))$. Therefore $K$ is a module of finite length and $\reg(K) \le \reg(G(\FilM))$.
 The   above exact sequence gives
 $$\begin{array}{ll}
 \reg (F_\qq(\FilM)/ x^*F_\qq(\FilM) ) & = \max\{\reg(K); \ \reg(F_\qq(\FilM/xM))\}\\
&  \le \max\{ \reg(G(\FilM)); \ \reg(F_\qq(\FilM/xM))\}.
\end{array}$$
  By \cite[Proposition 20.20]{E}, 
$$\reg (F_\qq(\FilM))  = \max\{ a_0(F_\qq(\FilM)) ,\ \reg(F_\qq(\FilM)/x^*F_\qq(\FilM))\}.$$
Hence
\begin{equation} \label{EF3}
\reg (F_\qq(\FilM) )  \le \max\{a_0(F_\qq(\FilM));\  \reg(G(\FilM)); \ \reg(F_\qq(\FilM/xM))\}.
\end{equation}
 Note that  $r(\FilM/xM) \le r$  and by Definition \ref{A1}(ii) $D(I,M/xM) \le D$. 
Using Theorem \ref{A3}, the inequality (\ref{EF3}) and Lemma \ref{F1} we get
 $$\begin{array}{ll}
 \reg (F_\qq(\FilM)) & \le  \max\{ (D+r+2)^2 - 2; \  (D+r+1)^2 - 1; \ 2D(D+r )+r -1 \}\\
& <  (D+r+2)^2 + D^2 - 3\end{array}$$
if $d=2$, 
$$ \begin{array}{ll}
\reg (F_\qq(\FilM))  & \le  \max\{ (D+r+2)^5 - 3; \  (D+r+1)^5 - 3; (D+r+2)^2 + D^2 -3 \} \\
& = (D+r+2)^5 - 3\end{array}$$
if $d=3$, and
$$ \begin{array}{ll}
\reg (F_\qq(\FilM)) & \le  \max\{ (D+r+2)^{3(d-1)!-1}- d; \  (D+r+1)^{3(d-1)!-1} - d; \\
& \hskip2cm (D+r+2)^{3(d-2)!-1}- d+1 \}\\
& = (D+r+2)^{3(d-1)!-1}- d\end{array}$$
for all $d\ge 4$.
 \end{pf}
 
  As an immediate consequence of the above theorem we get the following bound for the Castelnuovo-Mumford regularity of the classical fiber cone of an $\mm$-primary ideal.
 
 \begin{Corollary} \label{F4} Let 
 Let  $I$  be an $\mm$-primary ideal of $d$-dimensional local ring $A$. Then
\begin{itemize}
\item[(i)] $\reg(F_\mm (I)) \le 2 D(I,A)^2    -1$ if $d=1$;
\item[(ii)] $\reg(F_\mm(I)) \le 2 D(I,A)^2 +4 D(I,A) + 1$ if $d= 2$;
\item[(iii)] $\reg(F_\mm(I)) \le (D(I,A) +2)^{3(d-1)! - 1} - d$ if $d\ge 3$.
\end{itemize}
 \end{Corollary}
 
 In the graded case we can apply the method in Section \ref{Graded} to bound the Castelnuovo-Mumford regularity of $F_\qq(\FilM)$. We formulate here only  one result.
 
\begin{Proposition} Assume that $A$ is a homomorphic image of a Gorenstein graded algebra,  $M$ is a finitely generated graded $A$-module with $\dim M = d \geq  1$,  $I\subseteq \qq$ are graded $\mm$-primary  of $A$, and  $\FilM = \{ \Mc_n\}_{n\geq 0}$ is a good $I$-filtration of graded submodules of $M$  such that $\Mc_{n+1} \subseteq \qq \Mc_n$ for all $n\ge 0$.  Then
\begin{itemize}
\item[(i)] $\reg(F_\qq(\FilM)) \le  2 \ell(A/I) \hdeg (I,M)(\ell(A/I) \hdeg(I,M) + r(\FilM)) + r(\FilM) -1$ if $d=1$;
\item[(ii)] $\reg(F_\qq(\FilM)) \le (\ell(A/I)^2 \hdeg(I,M) +r(\FilM) +2)^2 + \ell(A/I)^4 \hdeg(I,M)^2 - 3$ if $d= 2$;
\item[(iii)] $\reg(F_\qq(\FilM)) \le (\ell(A/I)^d \hdeg(I,M) +r(\FilM)+2)^{3(d-1)! - 1} - d$ if $d\ge 3$.
\end{itemize}
 \end{Proposition}
 
 \vskip0.5cm

\noindent {\bf Acknowledgment}:  The authors would like to thank
the  referee for his/her valuable comments.

\end{document}